\documentclass[a4 paper,11pt,bezier]{article}
\usepackage[top=3cm,bottom=3cm,left=2.5cm,right=2.5cm]{geometry}
\usepackage{amsmath}
\usepackage{amsfonts,amsthm,amssymb}
\usepackage{amsfonts}
\usepackage{graphics,subfigure,caption2}
\usepackage{tikz, color}
\usepackage{cite}
\usepackage{graphics,epstopdf}
\usepackage{latexsym,bm}
\usepackage{amsfonts,amsthm,amssymb,mathrsfs,bbding}
\usepackage{CJK}
\usepackage{indentfirst}
\newtheorem{lem}{Lemma}[section]
\newtheorem{thm}[lem]{Theorem}

\theoremstyle{plain}

\begin{document}

\title{Graphs with a given conditional diameter that maximize the Wiener
	index\footnote{The research is supported by National Natural Science Foundation of China (12261086).}}
\author{Junfeng An, Yingzhi Tian\footnote{Corresponding author. E-mail: tianyzhxj@163.com (Y. Tian).}\\
{\small College of Mathematics and System Sciences, Xinjiang
University, Urumqi, Xinjiang 830046, China}}
\date{}
\maketitle
		
\maketitle {\flushleft\bf Abstract:} The Wiener index $W(G)$ of a graph $G$ is one of the most well-known topological indices, which is defined as the sum of distances between all pairs of vertices of $G$. The diameter $D(G)$ of $G$ is the maximum distance between all pairs of vertices of $G$; the conditional diameter $D(G;s)$ is the maximum distance between all pairs of vertex subsets with cardinality $s$ of $G$. When $s=1$, the conditional diameter $D(G;s)$ is just the diameter $D(G)$. The authors in \cite{QS} characterized the graphs with the maximum Wiener index among all graphs with diameter $D(G)=n-c$, where $1\le c\le 4$.
In this paper, we will  characterize the graphs with the maximum Wiener index among all graphs with conditional diameter $D(G;s)=n-2s-c$ ($-1\leq c\leq 1$), which extends partial results in \cite{QS}.

\maketitle {\flushleft\textit{\bf Keywords}:} Wiener index;  Diameter;  Conditional diameter
		
\section{Introduction}

Let $G$ be a simple  graph with vertex set $V(G)$ and edge set $E(G)$. The order and the size of $G$ are $n=|V(G)|$ and $m=|E(G)|$, respectively. The distance between two vertices $u$ and $v$, denoted by $d(u,v)=d_G(u,v)$, is the length of the shortest path connecting $u$ and $v$ in $G$. There are plenty of distance-based topological indices, which are widely used in mathematical chemistry in order to describe and predict the properties of chemical compounds. One of the most well-known topological indices is the Wiener index, which was introduced in 1947 by Wiener \cite{Wiener}. The 
Wiener index $W(G)$ of a graph $G$  is defined as the sum of distances between all (unordered) pairs of vertices of $G$, that is, $$W(G)=\sum _{\{u,v\}\subseteq V(G)} d(u,v).$$ 
Mathematical properties and applications of Wiener index are extensively studied, see \cite{QC, KC, AA, YL, Knor, Knor2, Cambie2, Cambie3, AV, DS, KX, XD} for references.

The diameter $D(G)$ of $G$ is the maximum distance between all pairs of vertices in $V(G)$, that is,  $D(G)=max_{u,v\in V(G)}d(u,v)$. For two  nonempty vertex subsets $V_1$ and $V_2$, the distance between $V_1$ and $V_2$, denoted by $d(V_1, V_2)=d_G(V_1, V_2)$, is the minimum of
the distances $d(x, y)$ among all $x \in V_1$ and $y \in V_2$.  Given a graphical property $\mathcal{P}$ satisfied by at least one pair $(V_1, V_2)$ of
nonempty subsets of $V(G)$, the conditional diameter $D_\mathcal{P} (G)$
of $G$  is
$$D_\mathcal{P}(G)=max\{d(V_1, V_2):\emptyset \ne V_1, V_2 \subseteq V(G),(V_1, V_2) \ satisfies \ \mathcal  {P}\}.$$
Note that $D_\mathcal{P} (G) = 0$ holds if and only if $V_1$ and $V_2$ overlap
for every $(V_1, V_2) \subseteq V(G) \times V(G)$ that satisfies $\mathcal{P}$. Conditional diameter measures the maximum distance between
subgraphs satisfying a given property. So their consideration
could be of some interest if in some applications we need to
control the communication delays between the network clusters modeled by such subgraphs.

The first choice of such a graphical property $\mathcal{P}_1$ is defined as follows:  ($V_1$, $V_2$) satisfies $\mathcal{P}_1$ if and only if $\left|V_{{1}}\right| = \left|V_{{2}}\right| =s$,  where $s$ is a positive integer. In this case, the conditional diameter  is denoted by $D(G;s)$, which is defined as 
$$ D(G;s) =  max\{d(V_{{1}} , V_{{2}}) : V_{{1}},V_{{2}} \subseteq V(G),\left|V_{{1}} \right| = \left|V_{{2}}\right| = {s}\}.$$
Clearly, $D(G; 1)$ is the standard diameter $D(G)$ of $G$. Thus $D(G;s)$ can be seen as a generalization of diameter $D(G)$. When $\left|V(G)\right|<2s$, then $D(G;s) = 0$. Moreover, when $ \left|V(G)\right|\ge 2s$, it is easy to see that the inequality $D(G;s) \le n-2s+1$ holds.

Although the Wiener index has been extensively studied, there are still some unsolved interesting questions. For example, Plesník \cite{Plesnik} asked an open problem ``What is the maximum average distance among graphs of order $n$ and diameter $d$?"; DeLaVi\~{n}a and Waller \cite{DeLaVina} conjectured that $W(G)\leq W(C_{2d+1})$ for any graph $G$ with  diameter $d\geq3$ and order $2d+1$, where $C_{2d+1}$ is the cycle of length $2d+1$. Some results related to the Wiener indices of graphs with given diameter can be seen in \cite{HL, SM, Wang}. Particularly, Cambie \cite{Cambie1} gave an asymptotic solution to the open problem of Plesník;
Sun et al. \cite{QS} characterized the graphs with the maximum Wiener index among all graphs with diameter $D(G)=n-c$, where $1\le c\le 4$. 

Motivated by the results above, we will investigate the maximum Wiener index among all graphs with given conditional diameter in this paper. Specifically,  we will characterize the graphs with the maximum Wiener index among all graphs with conditional diameter $D(G;s)=n-2s-c$, where $-1\leq c\leq 1$. Some lemmas will be given in the next section. Main results will be presented in the last section.

\section{Preliminaries}

The graphs considered in this paper are simple and undirected. For undefined notation and terminologies, we follow \cite{Bondy}. For a graph $G$, we denote by  $G-u$ and $G-uv$ the graphs obtained from $G$ by deleting the vertex $u\in V(G)$ and the edge $uv \in E(G)$, respectively. Similarly, $G + xy$ is a graph obtained from $G$ by adding an edge $xy \notin E(G)$. The induced subgraph $G[U]$ for a vertex subset $U\subseteq V(G)$ is $G-V(G)\setminus U$. The neighborhood of $u$ in $G$ is $N_G(u) = \{v | uv \in E(G)\}$. The degree $d_G(u)$ of $u$ in $G$ is $|N_G(v)|$. If $d_G(u)=1$, then $u$  is called a pendent vertex of $G$.
Denote by $P_n$ and $C_n$  the path and the cycle on $n$ vertices, respectively.

The sum of distances between $u$ and all other vertices of $G$ is $D_G(u) = \sum_{v\in V(G)}d(u,v)$.

\begin{lem}(\cite {AA}) 
	Let $G$ be a graph of order $n$, $v$ a pendent vertex of $G$ and $u$ the vertex adjacent to $v$. Then
	$W(G) = W(G-v)+D_{G-v}(u)+n-1$.
\end{lem}

\begin{lem} (\cite{HL})
	Let $P_n=v_1v_2\ldots v_n$ be a path on $n$ vertices, then $D_{P_n}(v_j)>D_{P_n}(v_{k})$  for $1\le j<k\le \lfloor n/2\rfloor$.
\end{lem}

	\begin{lem}(\cite{HL})
	Let $G$ be a non-trivial connected graph  on n vertices and let  $v\in V(G)$. Suppose that two paths $P=vv_1v_2\cdots v_k$, $Q=vu_1u_2\cdots u_l$ of lengths $k$, $l$ are attached to $G$ by their end vertices at $v$, respectively, to form $G_{k,l}$. As shown in Figure 1. If $l\ge k\ge 1$, then $W(G_{k,l}) < W(G_{k-1,l+1})$. 
\end{lem}

	\begin{figure}[htbp]
		\centering
		\includegraphics[width=10cm]{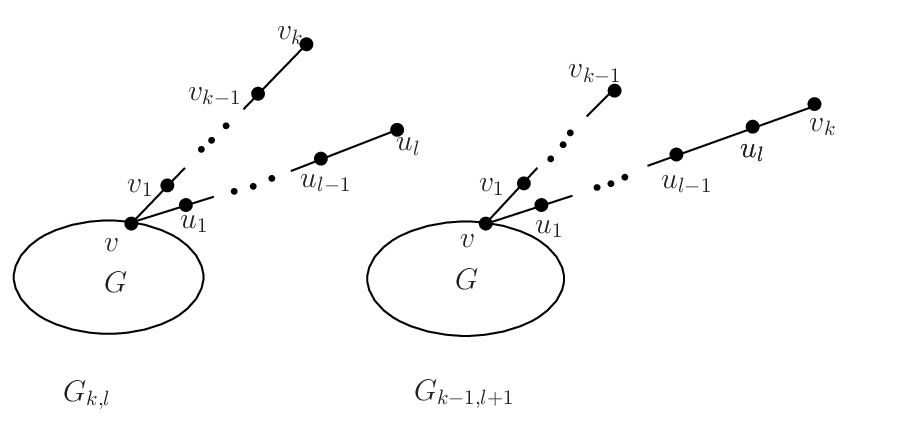}\\
		\caption{ $G_{k,l}$ and  $G_{k-1,l+1}$}
		\label{1}
	\end{figure}

\section{Main results}

It was proved in \cite{Knor} that $W(P_n)$ is maximum among all trees on $n$ vertices. Since removing 
of an edge from a connected graph results in increased  Wiener index, it is observed that Wiener index of a connected graph is less than or equal to the
Wiener index of its spanning tree. Thus $W(P_n)$ is maximum among all connected graph  on $n$ vertices. 
By $D(P_n;s)=n-2s+1$ ($n\ge 2s$), we have the following theorem.

\begin{thm}
	Let $G$ be a connected graph on $n$ vertices and $D(G;s)=n-2s+1$, where $s$ is a positive integer and  $n\ge 2s$. Then $W(G)\le W(P_n)$, and equality holds if and only if $G\cong P_n$.
\end{thm}

Let $T^{i}_n$ be the tree on $n$ vertices obtained from $P_{n-1}=x_1x_2\cdots x_{n-1}$ by attaching a pendent vertex  to $x_i$. See  Figure 2 for an illustration.

\begin{figure}[htbp]
	\centering
	\includegraphics[width=10cm]{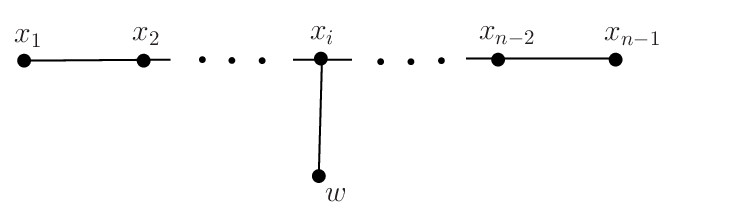}\\
	\caption{Tree $T^i_n$}
	\label{2}
\end{figure}

\begin{thm}
Let $G$ be a connected graph on $n$ vertices and $D(G;s)=n-2s$, where $s$ is a positive integer and  $n\ge 2s+3$.  Then $W(G)\le W(T^{s+1}_n)$, and equality holds if and only if $G\cong T^{s+1}_n$.
\end{thm}

\noindent{\bf Proof.} Let $d(L, R)=n-2s$, where $L=\{x_1,\ldots x_{s}\}$ and $R=\{x_{n-s}, \ldots, x_{n-1}\}$.
Assume $P=x_sx_{s+1}\cdots x_{n-s}$ is a  path of length $n-2s$ connecting $L$ and $R$ in $G$. Denote $M=\{x_{s+1},\ldots, x_{n-s-1}\}$ and $W=V(G)\setminus (L\cup M\cup R)=\{w\}$.

\noindent{\bf Claim.}  We can choose $L$ and $R$ such that $w$ is adjacent to vertices in $M$.

By $n\ge 2s+3$, $w$ can not adjacent to both vertices in $L$ and $R$. If $w$ is only adjacent to vertices in $L$, then $x_{s+1}$   must be adjacent to another vertex  $w' \in L$ other than $x_{s}$. Otherwise, the distance between $L-x_{s}+w$ and $R$ would be $n-2s+1$, a contradiction.
Thus we  replace $L$ by $L'$ and $W$ by $W'$,  where $ L'=L-w'+w$ and $W'=\{w'\}$. The case $w$ is only adjacent to vertices in $R$ can be analyzed similarly.   So the claim  holds.

\noindent {\bf Case 1.} $w$ is only adjacent to the vertices in $M$.

If $w$ is adjacent to more than one vertices in $M$, then delete all but one edges incident with $w$. Note that this operation does not change the conditional diameter and increases the Wiener index. So we assume  that $w$ is a pendent vertex. Without loss of generality, assume $w$ is adjacent to $x_i$, where  $s+1\le i\le n-s-1$.
 
Consider the induced subgraph $G[L\cup \{x_{s+1}\}]$. First, we transform  it to a tree by removing edges. Removing edges in this way does not change the conditional diameter and increases the Wiener index. Then, we transform it to a path as follows: we take one of the longest paths from $\{x_{s+1}\}$ and gradually enlarge it to an even longer path by appending the rest of the vertices in $L\cup \{x_{s+1}\}$  to the current endvertex on the other side of this path, one after another. By Lemma 2.3, each such transformation increases  the Wiener index and retains the conditional diameter. Similarly, we can transform $G[R\cup \{x_{n-s-1}\}]$ to a path with one endvertex $\{x_{n-s-1}\}$.

Now the graph $G$ is change to the graph isomorphic to  $T^{i}_n$, where $s+1\le i\le n-s-1$. Let $T^{s+1}_n= T^{i}_n-x_iw+x_{s+1}w$. Since $T^{s+1}_n-w\cong T^{i}_n-w\cong P_{n-1}$. By Lemmas 2.1 and  2.2, we have $W(T^{i}_n)\leq W(T^{s+1}_n)$, and equality holds if and only if $i=s+1$. Thus $W(G)\leq W(T^{s+1}_n)$, and equality holds if and only if $G\cong W(T^{s+1}_n)$.
 
\noindent {\bf Case 2.} $w$ is adjacent to vertices in both $L$ ($R$) and $M$.

We only need to consider that $w$ is adjacent to vertices in both $L$ and $M$. Since $P=x_sx_{s+1}\cdots x_{n-s}$ is a shortest  path  connecting $L$ and $R$, we obtain that  $N_G(w)\cap \{x_{s+1},\ldots,x_{n-s}\}\subseteq\{x_{s+1},x_{s+2}\}$. Let $x'=x_{s+2}$ if $x_{s+2}\in N_G(w)$ and  let $x'=x_{s+1}$ otherwise.

If $x'=x_{s+2}$, then consider the induced subgraph $G[L\cup \{x_{s+1},x_{s+2},w\}]$. First, we change it to a tree by removing some edges in  $E(G[R\cup \{x_{s+1},x_{s+2},w\}])\setminus\{x_{s+1}x_{s+2},x_{s+2}w\}$. Then, we transform it to a path such that $x_{s+2}$ is adjacent to one endvertex of this path as follows:  we take one of the longest paths from $\{x_{s+2}\}$ and gradually enlarge it to an even longer path by appending the rest of the vertices in $L$  to the current endvertex on the other side of this path, one after another. Note that $x_{s+2}$ is still adjacent to vertices  $x_{s+1}$ and $w$, and one of $x_{s+1}$ and $w$ must be an endvertex of this path.  Now we change $G$ to a graph isomorphic to $T^{s+2}_n$. By Case 1, we get $W(G)\leq W(T^{s+1}_n)$.

If $x'=x_{s+1}$, then by a similar argument as above, we can change $G$ to a graph isomorphic to $T^{s+1}_n$. Thus $W(G)\leq W(T^{s+1}_n)$.
  
From the arguments above, we obtain that $W(G)\le W(T^{s+1}_n)$, and the equality holds if and only if $G\cong T^{s+1}_n$.
$\hfill \square $

Let $T^{i,j}_n$ be a tree on $n$ vertices obtained from   $P_{n-2}=x_1x_2\cdots x_{n-2}$ by attaching two pendent vertices  to  $x_i$ and $x_j$, respectively. See  Figure 3 for an illustration.

\begin{figure}[htbp]
	\centering
	\includegraphics[width=12cm]{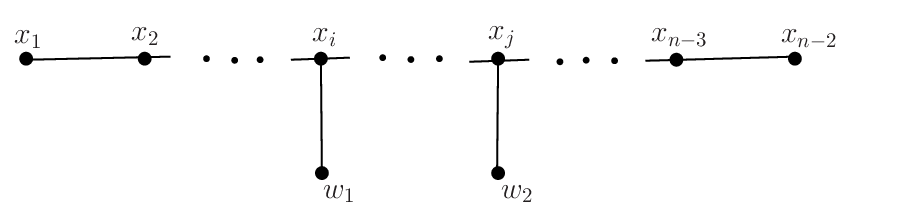}\\
	\caption{Tree $T^{i,j}_n$}
	\label{3}
\end{figure}

Let $T^{i(2)}_n$ be a tree on $n$ vertices obtained from $P_{n-2}=x_1x_2\cdots x_{n-2}$ by attaching the endvertex of a path of order 2 to $x_i$. See  Figure 4 for an illustration.

\begin{figure}[htbp]
	\centering
	\includegraphics[width=10cm]{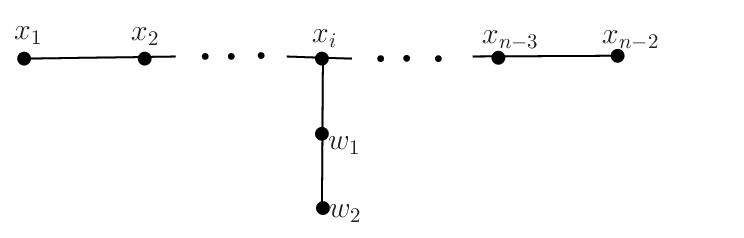}\\
	\caption{Tree $T^{i(2)}_n$}
	\label{4}
\end{figure}

\begin{thm}
	Let $G$ be a connected graph on $n$ vertices and $D(G;s)=n-2s-1$, where $s$ is a positive integer and  $n\ge 2s+5$. Then $W(G)\le W(T^{s+1,n-s-2}_n)$, and equality holds if and only if $G\cong T^{s+1,n-s-2}_n$.
\end{thm}

\noindent{\bf Proof.} 
Let $d(L, R)=n-2s-1$, where $L=\{x_1,\ldots x_{s}\}$ and $R=\{x_{n-s-1}, \ldots, x_{n-2}\}$.
Assume $P=x_sx_{s+1}\cdots x_{n-s-1}$ is a  path of length $n-2s-1$ connecting $L$ and $R$. Denote $M=\{x_{s+1},\ldots, x_{n-s-2}\}$ and $W=V(G)\setminus (L\cup M\cup R)=\{w_1,w_2\}$.

\noindent {\bf Case 1.} Neither $w_1$ nor $w_2$ is adjacent to vertices in $L\cup R$. 
 
\noindent {\bf Subcase 1.1.}  $w_1w_2\notin E(G)$.

Note that $N_G(w_i)\subseteq \{x_{s+1},\cdots, x_{n-s-2}\}$ for $i=1,2$.
If $w_i$ is adjacent to more than one vertices in $M$, then delete all but one edges incident with $w_i$, where $i\in\{1,2\}$. Note that this operation does not change the conditional diameter and increases the Wiener index. So we assume  that $w_i$ is pendent vertex for $i=1,2$. Without loss of generality, assume  that $w_1$is attached to $x_a$ and $w_2$ is attached to $x_b$, where $s+1\le a\le b\le n-s-2$.

By a similar argument as the proof of Case 1 in the Theorem 3.2, we transform $G[L\cup \{x_{s+1}\}]$ to a path with one endvertex $x_{s+1}$, and transform $G[R\cup \{x_{n-s-2}\}]$ to a path with one endvertex $x_{n-s-2}$.
That is, we change $G$ to a graph isomorphic to $T^{a,b}_n$ , where $s+1\le a\le b\le n-s-2$.

Let $T^{a,n-s-2}_{n}= T^{a,b}_n-x_bw_2+x_{n-s-2}w_2$, $k_{w_1}=d(x_{s+1} ,x_a)$ and $ k_{w_2}=d(x_{n-s-2} ,x_b)$. Since $T^{a,n-s-2}_{n}-w_2\cong  T^{a,b}_n-w_2$, we have
\begin{equation}
	\begin{split}
		\nonumber &W(T^{a,n-s-2}_{n})-W(T^{a,b}_n)\\
		&=1+2+\cdots +(n-s-2)+2+\cdots +(s+1)+(n-s-a)\\
		&-1-2-\cdots -(n-s-2-k_{w_2})-2-\cdots -(s+1+k_{w_2})-(n-s-a-k_{w_2})\\
		&=\sum_{1\le i\le k_{w_2}}(n-s-2-k_{w_2}+i)-\sum_{1\le i\le k_{w_2}}(s+1+i)+k_{w_2}.
	\end{split}
\end{equation}
Since  $n-2s-3-k_{w_2}\geq0$, we have
$ W(T^{i,n-m-2}_{n})-W(T^{i,j}_n)\ge 0$, and equality holds if and only if $k_{w_2}=0$.

Let $T^{s+1,n-s-2}_{n}=T^{a,n-s-2}_{n}-x_aw_1+x_{s+1}w_1$. Since $T^{s+1,n-s-2}_{n}-w_1\cong T^{a,n-s-2}_{n}-w_1$, we have
\begin{equation}
	\begin{split}
		\nonumber &W(T^{s+1,n-s-2}_{n})-W(T^{a,n-s-2}_{n})\\
		&=1+2+\cdots +(s+1)+2+\cdots +(n-s-2)+(n-2s-1)\\
		&-1-2-\cdots -(s+1+k_{w_1})-2-\cdots -(n-s-2-k_{w_1})-(n-2s-1-k_{w_1})\\
		&=-\sum_{1\le i\le k_{w_1}}(s+1+i)+\sum_{1\le i\le k_{w_1}}(n-s-2-k_{w_1}+i)+k_{w_1}.
	\end{split}
\end{equation}
By $n-2s-3-k_{w_2}\geq0$, we have 
$W(T^{s+1,n-s-2}_{n})-W(T^{a,n-s-2}_{n})\ge 0$, and equality holds if and only if $k_{w_1}=0$.

In this subcase, we conclude that $W(G)\le W(T^{s+1,n-s-2}_n)$, and equality holds if and only if $G\cong T^{s+1,n-s-2}_n$.

\noindent {\bf Subcase 1.2.} $w_1w_2\in E(G)$.

If both $w_1$ and $w_2$ have neighbors in $M$, then by deleting edge $w_1w_2$, we reduce this situation to Subcase 1.1.
So we assume  that only one of the vertices in $W$, say $w_1$, has neighbors in $M$ and the other vertex $w_2$ 
is a pendent vertex adjacent to $w_1$.  If  $w_1$ has more than one neighbors in $W$, then deleting all but one edges incident with $w_1$. Here the remaining edge satisfies the property that the end other than $w_1$ is farthest to the vertex set $\{x_{s+1},x_{n-s-2}\}$. Assume, without loss of generality, that $w_1$ is attached to $x_i$, where $s+2\le i\le n-s-3$.

By a similar argument as the proof of Case 1 in the Theorem 3.2, we transform $G[L\cup \{x_{s+1}\}]$ to a path with one endvertex $x_{s+1}$, and transform $G[R\cup \{x_{n-s-2}\}]$ to a path with one endvertex $x_{n-s-2}$.
That is, we change $G$ to a graph isomorphic to  $T^{i(2)}_n$,  where $s+2\le i\le n-s-3$.

Let $T^{(s+2)(2)}_{n}= T^{i(2)}_n-x_iw_1+x_{s+2}w_1$ and $k^{'}_{w_1}=d(x_{s+2} ,x_i)$.
Since $T^{(s+2)(2)}_{n}-w_1-w_2\cong  T^{i(2)}_n-w_1-w_2$, we have 
\begin{equation}
	\begin{split}
		\nonumber W(T^{(s+1)(2)}_{n})-W(T^{i(2)}_n)
		&=1+2+\cdots +(n-s-1)+3+\cdots +(s+2)\\
        &+1+2+\cdots +(n-s-2)+2+\cdots +(s+1)\\
		&-1-2-\cdots -(n-s-1-k^{'}_{w_1})-3-\cdots -(s+2+k^{'}_{w_1})\\
        &-1-2-\cdots -(n-s-2-k^{'}_{w_1})-2-\cdots -(s+1+k^{'}_{w_1})\\
		&=\sum_{1\le j\le k^{'}_{w_1}}(n-s-1-k^{'}_{w_1}+j)-\sum_{1\le j\le k^{'}_{w_1}}(s+2+j)\\
        &+\sum_{1\le j\le k^{'}_{w_1}}(n-s-2-k^{'}_{w_1}+j)-\sum_{1\le j\le k^{'}_{w_1}}(s+1+j).
	\end{split}
\end{equation}
Since $n-2s-3-k^{'}_{w_1}\geq 0$, we have
$W(T^{(s+2)(2)}_{n})-W(T^{i(2)}_n)\ge 0$, and equality holds if and only if $k^{'}_{w_1}=0$.

In this subcace, we conclude that $W(G)\leq W(T^{(s+2)(2)}_{n})$, and equality holds if and only if $G\cong T^{(s+2)(2)}_{n}$.

\noindent {\bf Case 2.} Either $w_1$ or $w_2$ are adjacent to vertices in $L\cup R$. 

If $w_i$ is only adjacent to vertices in $L\cup R$, then we can choose $L$ and $R$ such that $w_i$ is adjacent to some vertices in $M$ for $i=1,2$. Owing to Case 1, we only need to consider three subcases in the following.

\noindent {\bf Subcase 2.1.} Only one of $w_1$ or $w_2$, say $w_1$ is adjacent to vertices in $L\cup R$. 

We only consider $w_1$ is adjacent to  vertices in both $L$ and $M$, and $w_2$ is not adjacent to any vertices in $L\cup R$.

Since $P=x_sx_{s+1}\cdots x_{n-s-1}$ is a shortest  path  connecting $L$ and $R$, we obtain that  $N_G(w_1)\cap \{x_{s+1},\ldots,x_{n-s}\}\subseteq\{x_{s+1},x_{s+2}\}$. Let $x'=x_{s+2}$ if $x_{s+2}\in N_G(w_1)$ and  let $x'=x_{s+1}$ otherwise. If $x'=x_{s+2}$, then by a similar argument as the proof of Case 2 in Theorem 3.2, we can change $G[L\cup \{x_{s+1},x_{s+2},w\}]$ to a path such that $x_{s+2}$ is still adjacent to vertices  $x_{s+1}$ and $w$, and one of $x_{s+1}$ and $w_1$ is an endvertex of this path. If $x'=x_{s+1}$, then by a similar argument, we can change $G[L\cup \{x_{s+1},w\}]$ to a path such that $x_{s+1}$ is still adjacent to vertices  $x_{s}$ and $w_1$, and one of $x_{s}$ and $w_1$ is an endvertex of this path.

Suppose $w_2$ is adjacent to some vertices in $M$, then deleting all edges incident with $w_2$ but one edge joining $w_2$ to a vertex in $M$. Then $G$ is changed to a graph isomorphic to $T_n^{i,j}$. Suppose $w_2$ is only adjacent to $w_1$, then $G$ is changed to a graph isomorphic to $T_n^{i(2)}$.

\noindent {\bf Subcase 2.2.} Both $w_1$ and $w_2$ are adjacent to vertices in $L$ ($R$). 

We only consider $w_i$ is adjacent to vertices in  both $L$ and $M$ for $i=1,2$.

Since $P=x_sx_{s+1}\cdots x_{n-s-1}$ is a shortest  path  connecting $L$ and $R$, we obtain that  $N_G(w_i)\cap \{x_{s+1},\ldots,x_{n-s}\}\subseteq\{x_{s+1},x_{s+2}\}$ for $i=1,2$. Let $x'=x_{s+2}$ if $x_{s+2}\in N_G(w_1)$ and  let $x'=x_{s+1}$ otherwise. Let $x''=x_{s+2}$ if $x_{s+2}\in N_G(w_2)$ and  let $x''=x_{s+1}$ otherwise. Here we only give the proof when  $x'=x_{s+2}$ and $x''=x_{s+2}$. Other cases can be proved similarly.
We consider the induced subgraph $G[L\cup \{x_{s+1},x_{s+2},w_1,w_2\}]$. First, we change it to a tree by removing some edges in  $E(G[L\cup \{x_{m+1},x_{m+2},w_1,w_2\}])\setminus\{x_{s+1}x_{s+2},x_{s+2}w_1,x_{s+2}w_2\}$. Then, we transform it to a tree such that $x_{s+2}$ is adjacent to two pendent vertices as follows:  we take one of the longest paths from $x_{s+2}$ and gradually enlarge it to an even longer path by appending the rest of the vertices in $L$  to the current endvertex on the other side of this path, one after another. Note that $x_{s+2}$ is still adjacent to vertices  $x_{s+1}$, $w_1$ and $w_2$, and two of  $x_{s+1}$, $w_1$ and $w_2$ are pendent vertices adjacent to $x_{s+2}$.  Then $G$ is changed to a graph isomorphic to $T_n^{i,j}$.

\noindent {\bf Subcase 2.3.} One of $w_1$ or $w_2$, say $w_1$, is adjacent to vertices in $L$ and $w_2$ is adjacent to vertices in  $R$. 

Since $P=x_sx_{s+1}\cdots x_{n-s-1}$ is a shortest  path  connecting $L$ and $R$, we obtain that  $N_G(w_1)\cap \{x_{s+1},\ldots,x_{n-s}\}\subseteq\{x_{s+1},x_{s+2}\}$. Let $x'=x_{s+2}$ if $x_{s+2}\in N_G(w_1)$ and  let $x'=x_{s+1}$ otherwise. If $x'=x_{s+2}$, then by a similar argument as the proof of Case 2 in Theorem 3.2, we can change $G[L\cup \{x_{s+1},x_{s+2},w_1\}]$ to a path such that $x_{s+2}$ is still adjacent to vertices  $x_{s+1}$ and $w_1$, and one of $x_{s+1}$ and $w_1$ is an endvertex of this path. If $x'=x_{s+1}$, then by a similar argument, we can change $G[L\cup \{x_{s+1},w_1\}]$ to a path such that $x_{s+1}$ is still adjacent to vertices  $x_{s}$ and $w_1$, and one of $x_{s}$ and $w_1$ is an endvertex of this path.
Similarly, if $w_2x_{n-s-3}\in E(G)$, we can change $G[R\cup \{x_{n-s-3},x_{n-s-2},w_2\}]$ to a path that $x_{n-s-3}$ is still adjacent to vertices  $x_{n-s-2}$ and $w_2$, and one of $x_{n-s-2}$ and $w_2$ is an endvertex of this path. If $w_2x_{n-s-3}\notin E(G)$, we can change $G[R\cup \{x_{n-s-2},w_2\}]$ to a path that $x_{n-s-2}$ is still adjacent to vertices  $x_{n-s-1}$ and $w_2$, and one of $x_{n-s-1}$ and $w_2$ is an endvertex of this path. Thus $G$ is changed  to a graph isomorphic to $T_n^{i,j}$.

All cases lead to $W(G)\le W(T^{s+1,n-s-2}_n)$ or $W(G)\le W(T^{(s+2)(2)}_n)$. So we only need to compare $W(T^{s+1,n-s-2}_n)$ and $W(T^{(s+2)(2)}_n)$. Since 
$W(T^{s+1,n-s-2}_n)-W(T^{(s+2)(2)}_n)=\frac{1}{2}n^2-(s+\frac{3}{2})n+s^2+7s>0$, we obtain that 
$W(G)\le W(T^{s+1,n-s-2}_n)$, and equality holds if and only if $G\cong T^{s+1,n-s-2}_n$.
The proof is thus complete.
$\hfill \square $

\end{document}